\documentclass{article}
\usepackage{maa-monthly}
\usepackage[outdir=./]{epstopdf}

\begin{document}

\title{Expected Number of Vertices of a Hypercube Slice}
\markright{Expected Number of Vertices of a Hypercube Slice}
\author{Hunter Swan}

\maketitle

\begin{abstract}
Given a random k-dimensional cross-section of a hypercube, what is its expected number of vertices? We show that, for a suitable distribution of random slices, the answer is $2^k$, independent of the dimension of the hypercube.
\end{abstract}

To understand a high-dimensional object, a good starting point is to take a look
at lower dimension cross-sections which can be more easily visualized. For this to be
effective, though, we have to have some point of comparison---i.e. an idea of what such
cross-sections look like for ``typical'' high dimensional objects. This note considers
what cross sections ``look like'' for the case of a hypercube. In particular, we ask:
What is the expected number of vertices for a random $k$-dimensional slice of an $n$-dimensional
hypercube? We will show that, for a suitable distribution of slices, the
answer is $2^k$, regardless of the dimension $n$ of the hypercube.

Throughout our discussion, we fix our hypercube to be $C := [-1, 1]^n$. A $k$-dimensional slice (i.e. cross-section) of the cube is the intersection of $C$ with a
$k$-dimensional affine subspace of $\mathbb{R}^n$ (henceforth referred to as a $k$-flat), and is characterized by the following two properties.
\begin{itemize}
\item An orientation, specified by $k$ unit vectors $n_1,\dots, n_k \in S^{n-1}$.
\item A translation, specified by some $\tau \in \mathbb{R}^n$, which may be taken WLOG to be perpendicular
to the unit vectors $n_i$.
\end{itemize}
The slice is then given explicitly by $C\cap (\tau + span\{n_1,\dots, n_k\})$.

There are various different ways to choose a random slice (i.e. distributions), but one
of the most natural is the so-called distribution of \textit{isotropic random k-flats} \cite{SchneiderWeil}.  This
comes from choosing the unit vectors $n_i$ from $S^{n-1}$ uniformly and independently and
then choosing the translation $\tau$ uniformly from all those translations for which the
resulting $k$-flat still intersects C. The result of this note holds for this distribution, and
in fact holds much more generally.  We will show that, for any fixed orientation, by
choosing the translation $\tau$ uniformly from all translations for which the resulting flat
still intersects the hypercube, the expected number of vertices of the resulting slice is
$2^k$. As a corollary, for \textit{any} distribution of orientations, so long as we then choose the
translation uniformly, the expected number of vertices is $2^k$.

To put this in perspective, we remark that a \textit{particular} 2-slice of a $n$-dimensional
hypercube (for $n > 2$) can have anywhere from 3 to $2n$ vertices. An interpretation of
this result for the case $k = 2$ is thus that there is a very small chance of finding a slice
of a hypercube with very many vertices.

To compute the expected number of vertices, first fix the orientation of a $k$-slice, as
given by (linearly independent) unit vectors $n_1, \dots, n_k$. The translation $\tau$ may then be
taken to be a vector in the $(n - k)$-dimensional subspace normal to $n_1,\dots, n_k$, which
we denote by $N$. The range of values $\tau$ may take so that the slice still intersects the
hypercube $C$ is given by the projection $P_N (C)$ of $C$ on to $N$, where $P_N (\cdot)$ is the
projection operator on to $N$.

\begin{figure}[h]
	\begin{centering}
	\includegraphics[width=.5\columnwidth]{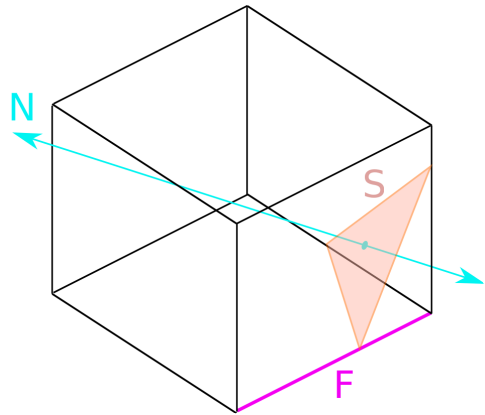}
	\caption[2-slice of a 3-cube]{A typical 2-dimensional slice $S$ of a 3-cube. The vertices of the slice result from intersections with
1-dimensional faces (edges) of the cube, such as the one labeled $F$. $N$ is the 1-dimensional subspace normal to $S$.
	}
	\label{slice}
	\end{centering}
\end{figure}

A vertex results from the intersection of a $k$-flat with an $(n - k)$-dimensional face
$F$ of $C$. For a given translation $\tau$, the resulting slice $(\tau + span\{n_1,\dots, n_k\}) \cap C$
intersects a particular $(n - k)$-dimensional face $F$ of $C$ if and only if $\tau \in P_N (F)$. Since $\tau$ is
assumed to be chosen uniformly within $P_N (C)$, the probability that the slice intersects
$F$ is just the ratio of the $(n - k)$-volume of $P_N (F)$ to that of $P_N (C)$. We will need a
formula for the volume of these regions, and to accomplish this we first introduce the
concept of a \textit{zonotope}.

Given arbitrary vectors $v_1, \dots, v_m \in \mathbb{R}^n$, the \textit{zonotope} $Z (v_1, \dots, v_m)$ is defined as
\begin{equation*}
Z (v_1,\dots, v_m) := \{\lambda_1 v_1 +\cdots + \lambda_m v_m | \lambda_1,\dots, \lambda_m \in [0, 1]\}.
\end{equation*}
Translations of such sets we also refer to as zonotopes. For a set of vectors $S =\{v_1, \dots, v_m\}$ we also use the notation $Z(S)$ to refer to the zonotope $Z(v_1, \dots, v_m)$.
When the vectors $v_i$ are linearly independent, the zonotope $Z(v_1, \dots, v_m)$ is just a parallelotope,
and conversely any parallelotope is a zonotope, but in general a zonotope
is a more complicated object (e.g. when $m > n$). We note that all of $C$, $F$, $P_N (C)$,
and $P_N (F)$ are zonotopes.  For $C$ and $F$ this is obvious, because they are both parallelotopes.
For the other two, we see from the linearity of $P_N (\cdot)$ that the projection of
a zonotope is again a zonotope, i.e.
\begin{align*}
P_N (Z(v_1, \dots, v_m)) & = \{\lambda_1 P_N(v_1) +\cdots+ \lambda_m P_N(v_m) | \lambda_1,\dots, \lambda_m \in [0, 1]\} \\
& = Z(P_N (v_1),\dots, P_N(v_m)).
\end{align*}

As an aside, this fact that a zonotope may be realized as the projection of a parallelotope
is true in general.  Any zonotope is the projection of some parallelotope (which is
not hard to see from the definition we gave above), and this can actually be taken as
a \textit{definition} of a zonotope\footnote{An informal interpretation of this definition is that ``a zonotope is the shadow cast by a box of arbitrary
orientation''.}.  From this perspective it is obvious why zonotopes arise in
the current problem, where we are concerned with projections of the parallelotopes $C$
and $F$ into the space $N$.

\begin{figure}[h]
	\begin{centering}
	\includegraphics[width=.5\columnwidth]{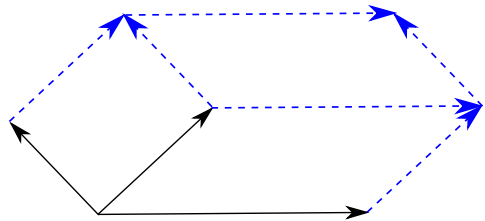}
	\caption[A zonotope]{A 2-dimensional zonotope and its partitioning into parallelotopes (parallelograms, in this case). The
solid black arrows are the vectors defining the zonotope, and the zonotope itself is the entire hexagonal region.
The smaller parallelogram regions outlined by the blue and black arrows are those that partition the zonotope.
Fun fact: this partitioning is not unique. Can you find another?
	}
	\label{zonotope}
	\end{centering}
\end{figure}

Zonotopes have a number of interesting combinatorial properties\footnote{One of the chief motivations for studying zonotopes comes from their combinatorial relation to arrangements
of hyperplanes \cite{McMullen}.}. Most importantly for our purposes, a zonotope $Z(v_1, \dots, v_m)$ which has dimension $d$ can be \textit{partitioned}
into a number of parallelotopes \cite{BeckRobins}, with one parallelotope for each maximal
linearly independent subset $S \subseteq \{v_1, \dots, v_m\}$. These partitioning parallelotopes have
the form $Z(S)$ for each such subset $S$. An example of this partitioning can be seen
in Figure \ref{zonotope}. This decomposition allows us to evaluate the $d$-volume of $Z(v_1, \dots, v_n)$
as the sum of the $d$-volumes of the partitioning parallelotopes $Z(S)$, leading to the
following formula for the volume of a zonotope \cite{Shephard,BeckRobins}:
\begin{equation}
V_d(Z(v_1,\dots, v_m)) = \sum_{\substack{S\subseteq \{v_1,\dots,v_m\} : \\ |S|=d}} V_d(Z(S)),
\end{equation}
where $V_d(\cdot)$ denotes $d$-dimensional volume. (Note that in this formula we sum over
$d$-subsets of $\{v_1, \dots, v_m\}$, whereas we expressed the partitioning of $Z(v_1, \dots, v_m)$ in
terms of maximal linearly independent subsets. In fact, we could just as well sum over
maximal linearly independent subsets in the volume formula, but the form we have
written will be more useful to us later. The equality of the two sums follows from the
fact that every maximal linearly independent subset has $d$ vectors, and any $d$-subset
which is not linearly independent yields a degenerate parallelotope with zero volume.)
The utility of this formula in general is that when $Z(S)$ is a parallelotope, its volume
can be found via the usual determinant formula\footnote{The fact that zonotopes admit such a concise analytical expression for their volume makes them very
special in the world of polytopes, where computing volumes is generally difficult \cite{DyerFrieze,BeckRobins}.}. In the present case, it turns out that
we will never actually need to evaluate a determinant, as all factors of $V_d(\cdot)$ will
eventually cancel.

The cube $C$ is, up to a translation, $Z(2e_1, \dots, 2e_n)$, where $e_i$ are the standard basis
vectors for $\mathbb{R}^n$. So $P_N(C) = Z(2 P_N(e_1), \dots, 2 P_N (e_n))$. Similarly, an $(n-k)$-dimensional face $F$ has the form (again up to translation) 
$Z(2 e_{i_1},\dots, 2 e_{i_{n-k}})$ for some $n-k$ standard basis vectors $e_{i_j}$, so $P_N(F) = Z(2 P_N (e_{i_1}),\dots, 2 P_N (e_{i_{n-k}}))$.
As mentioned above, the probability that face $F$ is intersected by a random flat is given by the ratio of the volumes
\begin{equation*}
\frac{V_{n-k}(P_N (F))}{V_{n-k}(P_N (C))} = \frac{V_{n-k}(Z(2 P_N (e_{i_1}),\dots, 2 P_N(e_{i_{n-k}})))}{V_{n-k}(Z(2 P_N (e_1), \dots, 2 P_N (e_n)))}.
\end{equation*}
The expected number of vertices, which we denote as \#, is the sum
\begin{equation*}
\# = \sum_F \frac{V_{n-k}(P_N (F))}{V_{n-k}(P_N (C))} = \frac{\sum_F V_{n-k}(P_N (F))}{V_{n-k}(P_N (C))}
\end{equation*}
of this probability over all $(n-k)$-dimensional faces $F$. Note that for any set of $n-k$
standard basis vectors $e_{i_1},\dots, e_{i_{n-k}}$ there are $2^k$ faces $F$ which are translations of the
zonotope $Z(2 e_{i_1},\dots, 2 e_{i_{n-k}})$; this is because such a face $F$ has $k$ fixed coordinates $x_i$, 
each of which is either $+1$ or $-1$. Since each such face has the same probability of
being intersected by a random flat, we can replace the sum over faces $F$ in the above
equation by a sum over $(n-k)$-subsets of $\{2 e_1, \dots, 2 e_n\}$:
\begin{equation*}
\# = \sum_{\substack{ S\subseteq \{2e_1,\dots,2e_n\}: \\ |S|=n-k}} \frac{2^k V_{n-k}(Z(P_N (S)))}{V_{n-k}(P_N(C))}.
\end{equation*}

The only thing that remains is to expand $V_{n-k}(P_N (C))$ using the volume formula
for zonotopes:
\begin{align*}
\# &= \left(\frac{ \sum\limits_{\substack{ S\subseteq \{2e_1,\dots,2e_n\}: \\ |S|=n-k}} 2^k V_{n-k}(Z(P_N (S)))}{V_{n-k}(Z(P_N (2 e_1),\dots, P_N (2 e_n)))} \right) \\
&= 2^k  \left(\frac{ \sum\limits_{\substack{ S\subseteq \{2e_1,\dots,2e_n\}: \\ |S|=n-k}} V_{n-k}(Z(P_N (S)))}{ \sum\limits_{\substack{ S\subseteq \{2e_1,\dots,2e_n\}: \\ |S|=n-k}} V_{n-k}(Z(P_N (S)))} \right) \\
&= 2^k.
\end{align*}

Which is our desired result.

Note that this argument generalizes immediately to the case where the hypercube $C$ is replaced by an arbitrary parallelotope.

\begin{acknowledgment}{Acknowledgment.}
This is the author's first paper, and he would like to use the occassion to acknowledge his greatest Sponsor in the words of Proverbs 3:5-6: \textit{Trust in the Lord with all your heart, and lean not on your own understanding.  In all your ways acknowledge him, and he will make your paths straight.}
\end{acknowledgment}

\end{document}